# Performance assessment of two active power filter control strategies in the presence of non-stationary currents

Vin Cent Tai, Phen Chiak See, Marta Molinas, and Kjetil Uhlen

Norwegian University of Science and Technology

This paper describes an active power filter (APF) control strategy, which eliminates harmonics and compensates reactive power in a three-phase four-wire power system supplying non-linear unbalanced loads in the presence of non-linear non-stationary currents. Empirical Mode Decomposition (EMD) technique developed as part of the Hilbert-Huang Transform (HHT) is used to singulate the harmonics and non-linear non stationary disturbances from the load currents. The control strategy for the APF is formulated by hybridizing the so called modified *p-q* theory with the EMD algorithm. A four-leg split-capacitor converter controlled by hysteresis band current controller is used as an APF. The results obtained are compared with those obtained with the conventional modified *p-q* theory, which does not possess current harmonics or distortions separation strategy, to validate its performance.

Last updated on 7-Jun-12

## 1   Introduction

Non-linear non-stationary currents (also known as non-periodic currents) are known to be problematic to power systems. Arc furnaces, cycloconverters, adjustable speed drives, as well as transient disturbances are the typical sources that generate non-linear non-stationary currents [Czarnecki, 2000; and Tolbert et al., 2003]. This also includes currents with sub-harmonics as well as super-harmonics [Tlusty et al., 2012]. The compensation with the conventional power theories becomes erratic due to the energy flow in the presence of non-linear non-stationary currents differs from that in the presence of only reactive, unbalanced, and harmonic currents [Czarnecki, 2000].

Czarnecki [Czarnecki, 2000] suggested that non-linear non-stationary currents should be compensated along with harmonics by the same device that reduces the reactive and unbalanced currents. He also pointed out that this can be achieved through the proper choice of measurement and digital signal processing procedure as well as the right control algorithm. For instance, Prony analysis is used in [Nakajima et al., 1988; and Qi et al., 2006] to formulate the APF control for eliminating the stationary and non-stationary harmonics; Generalized non-active power theory (GNP), a theory based on the Fryze's current theory has been formulated in [Peng and Tolbert, 2000] to cope with the such current components; and Wavelet analysis is employed to analyze the performance of the GNP [Tlusty et al., 2012].

This paper focuses on the formulation of the control strategy for an active power filter, which employs EMD technique and modified $p$-$q$ theory to achieve non-linear non-stationary disturbances and harmonics elimination in the source currents, as well as reactive power compensation. Unlike the conventional compensation strategy which does not employ distortion singulation technique, the proposed strategy is able to compensate not only harmonic currents and reactive power, but also non-linear non-stationary disturbances from the source currents.

The first section of this paper introduces the EMD and the modified $p$-$q$ theory used in the proposed control strategy. Then, the formulation of the control strategy is addressed. Subsequently, a test case with perfect three-phase voltage source supplying unbalanced load with non-linear non-stationary disturbances is presented. The results obtained with the proposed control strategy are then compared with the ones obtained with the conventional modified $p$-$q$ theory. Finally, this paper ends with conclusions.

## 2   Compensation methods

The authors propose to combine EMD technique with modified $p$-$q$ theory to compensate for reactive power, harmonics, and distortion caused by non-linear non-stationary currents.

## 2.1 Empirical mode decomposition

EMD was introduced by [Huang et al., 1998] to decompose non-linear non-stationary time series data into several zero-mean components known as Intrinsic Mode Functions (IMFs). In general, EMD works as follow, with input signal in time series, $x(t)$:

---
Algorithm 1: EMD

---

1. Identify all local extrema of $x(t)$.
2. Create the upper envelop by joining the maxima with the interpolated cubic spline.
3. Create the lower envelop for the local minima with the interpolated cubic spline.
4. Compute the mean value, $m_i$ from the upper and lower envelop.
5. Compute the $i^{th}$ sift component, $h_i = x(t) - m_i$.
6. Iterate steps 1 to 5 until $h_i$ possesses zero mean and free from any extrema. The signal obtained is call an IMF and is designated as $c_k$, after $k^{th}$ iteration.
7. Compute the residue, $r_k = x(t) - c_k$.

---

The number of extrema decreases after every iteration from one $r_k$ to next. The above steps are repeated until a predetermined stopping criterion is met. Standard deviation (SD) is used as the stopping criteria. The value of the SD is typically between 0.2 and 0.3 [Huang et al., 1998].

$$SD = \sum_{t=0}^{T} \left[\frac{|h_{i-1}(t) - h_i(t)|^2}{h_{i-1}^2(t)}\right] \tag{1}$$

In the following analysis, these IMFs are the current components (non-periodic, harmonics, sub-harmonics, super-harmonics etc.) that caused distortions in the power system, which needed to be eliminated.

## 2.2 Modified $p$-$q$ power theory

The $p$-$q$ theory (or instantaneous power theory) is the most widely used power theory to formulate the control of APF [Tlusty et al., 2012]. To compensate the zero sequence current and zero sequence power that exist in three-phase four-wire systems under unbalance-load condition, the modified $p$-$q$ theory [Togasawa et al., 1994; Akagi et al., 1999] is adapted to calculate the compensating currents. The modified $p$-$q$ theory is the power theory developed on the basic of $p$-$q$ theory, which takes into account the zero sequence components that occurs in a three-phase four-wire system under unbalanced load conditions. In $\alpha - \beta - 0$ coordinates, the line voltages and currents are:

$$\begin{bmatrix} v_\alpha \\ v_\beta \\ v_0 \end{bmatrix} = \sqrt{\frac{2}{3}} \begin{bmatrix} 1 & -1/2 & -1/2 \\ 0 & \sqrt{3}/2 & -\sqrt{3}/2 \\ 1/\sqrt{2} & 1/\sqrt{2} & 1/\sqrt{2} \end{bmatrix} \begin{bmatrix} v_R \\ v_S \\ v_T \end{bmatrix} \qquad (2)$$

$$\begin{bmatrix} i_\alpha \\ i_\beta \\ i_0 \end{bmatrix} = \sqrt{\frac{2}{3}} \begin{bmatrix} 1 & -1/2 & -1/2 \\ 0 & \sqrt{3}/2 & -\sqrt{3}/2 \\ 1/\sqrt{2} & 1/\sqrt{2} & 1/\sqrt{2} \end{bmatrix} \begin{bmatrix} i_R \\ i_S \\ i_T \end{bmatrix} \qquad (3)$$

The instantaneous real power, $p$ is defined as the dot product of voltage vector $\boldsymbol{v}$ and current vector, $\boldsymbol{i}$, which can be further decomposed into its average and oscilating components, $\bar{p}$ and $\tilde{p}$. On the other hand, the cross product of $\boldsymbol{i}$ and $\boldsymbol{v}$ gives the instantaneous imaginary power as a vector $\boldsymbol{q}$, which can be further decomposed into $\alpha - \beta - 0$ components as $q_\alpha$, $q_\beta$, and $q_0$. The magnitude of the instantaneous imaginary power is norm of the vector, $\boldsymbol{q}$.

$$p = i_\alpha v_\alpha + i_\beta v_\beta + i_0 v_0 = \bar{p} + \tilde{p} \qquad (4)$$

$$\boldsymbol{q} = \begin{bmatrix} q_\alpha \\ q_\beta \\ q_0 \end{bmatrix} = \begin{bmatrix} 0 & -v_0 & v_\beta \\ v_0 & 0 & -v_\alpha \\ -v_\beta & v_\alpha & 0 \end{bmatrix} \begin{bmatrix} i_\alpha \\ i_\beta \\ i_0 \end{bmatrix} \qquad (5)$$

## 3 Control of the APF with EMD

The objective of the compensation is to eliminate the harmonic and non-linear non-stationary currents, and reactive power. Figure 1 shows the control scheme for the APF. The proposed control method treats the distortion elimination and reactive power compensation separately.

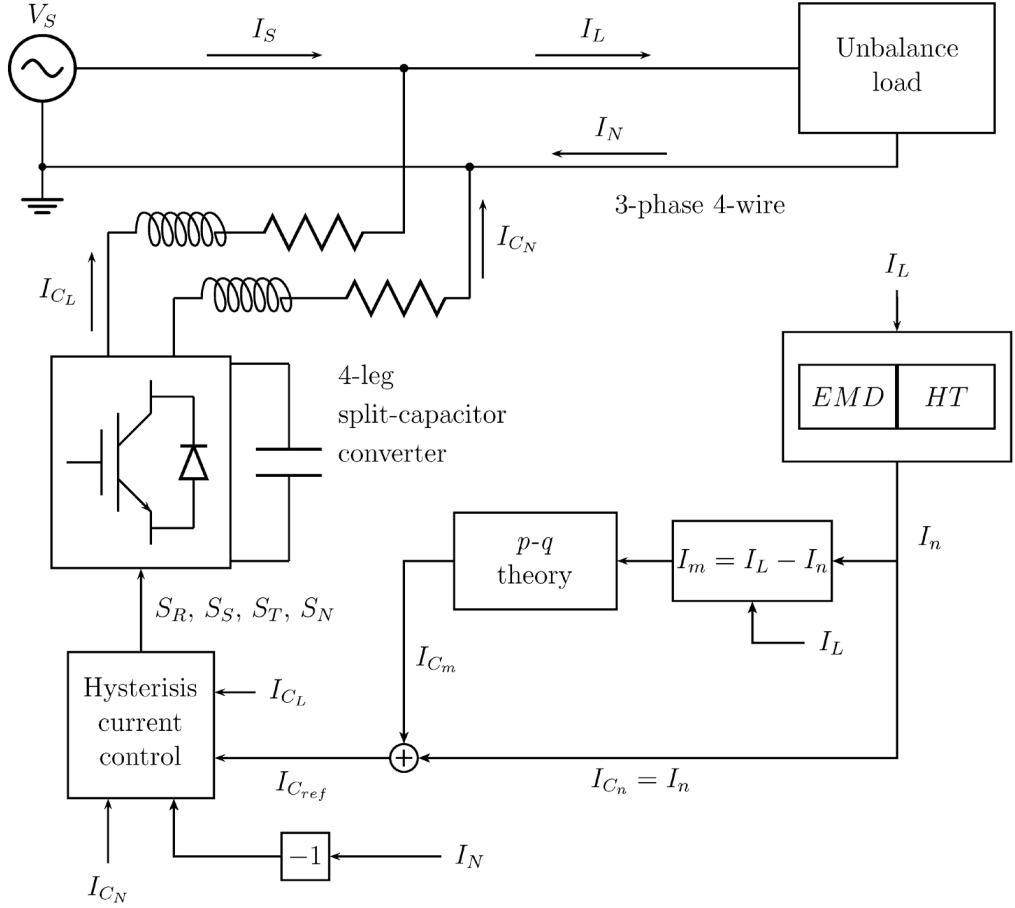

Figure 1. Proposed control scheme for the APF.

The filtering process starts with measurement of the line and neutral currents. The measured line current, $\boldsymbol{I_L} = [i_R,\ i_S,\ i_T]^T$ is then fed into the EMD block (see Figure 1) for signals decomposition.

Czarnecki [Czarnecki, 2000] suggested the load current can be decomposed into fundamental and residual components (in this paper, they are $\boldsymbol{I_m}$ and $\boldsymbol{I_n}$, respectively). For a system with balance sinusoidal source voltage, he pointed out $\boldsymbol{I_n}$ is composed of harmonics and non-linear non-stationary components, while $\boldsymbol{I_m}$ is the one that has an effect on the power factor reduction.

The residual current $I_n = [i_{nR}, i_{nS}, i_{nT}]^T$ is obtained through the EMD block. The fundamental current is calculated by subtracting the line current with the residual components [Czarnecki, 2000].

$$I_m = I_L - I_n \tag{6}$$

$I_m$ is then processed by using the modified $p$-$q$ theory to obtain the reference compensating currents, $I_{C_m}$ for reactive power compensation. The total reference compensating current is calculated by adding $I_{C_m}$ and $I_{C_n}$, where $I_{C_n} = I_n$:

$$I_{C_{ref}} = I_{C_m} - I_{C_n} \tag{7}$$

$I_{C_m}$ and $I_{C_n}$ are the compensating currents obtained by using $I_m$ and $I_n$, respectively. To obtain $I_{C_m}$, a second-order low-pass filter with cutoff frequency of 8Hz is used to extract the oscillating real power, $\tilde{p}$ that is needed to be compensated.

A four-leg split-capacitor converter is used to construct the shunt active filter. A hysteresis band current controller is used to generate the switching pulses to control the IGBTs of the converter. The switching logic is formulated as follows:

---
Algorithm 2: Hysteresis band current control

---
if $I_{C_L} > (I_{C_{ref}} + HB)$ then
    Upper switch is ON
    Lower switch is OFF
else
    Upper switch is OFF
    Lower switch is ON
end if

---

where, $HB$ is the hysteresis band, and $I_{C_L}$ is the compensating current generated by the APF. This algorithm allows $I_{C_L}$ to follow $I_{C_{ref}}$ within the preset band at all time.

The neutral current is compensated simply by generating the neutral current in the reverse direction. The results anticipated at the end of the compensation are distortion-free sinusoidal source currents with unity power factor, and zero reactive power from the power source.

## 4 Simulation results

The power system chosen for this study is shown in Figure 1. The voltage source is supplying the three-phase four-wire system 110V at 50Hz, and the resistance and inductance of the line ($R_L$ and $L_L$) are 0.07Ω and 45mH, respectively. The disturbance current is then injected into line 1, at time between 0.088s and 0.094s, to simulate the instability caused by the non-linear non-stationary signal. The source volage and the resultant load current waveforms are shown in Figure 2.

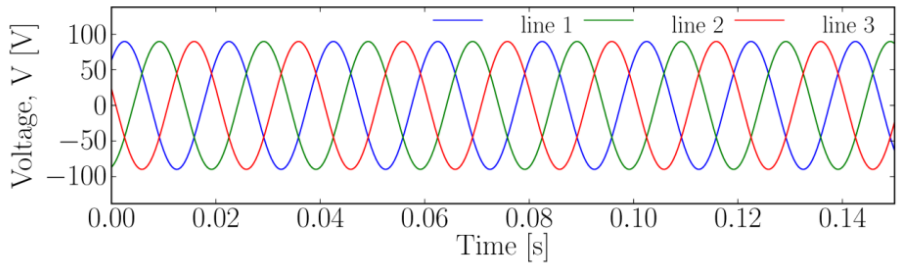

(a)

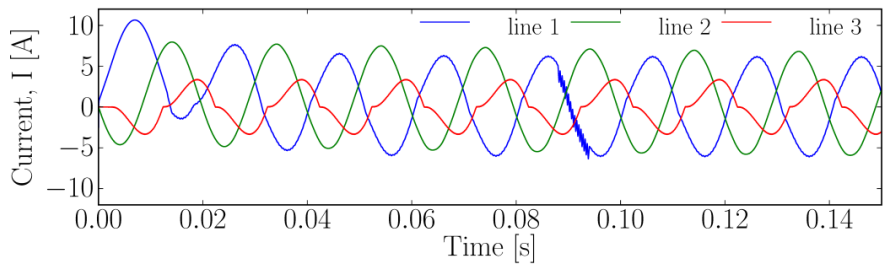

(b)

Figure 2. Simulation input: (a) source voltages; and (b) load currents.

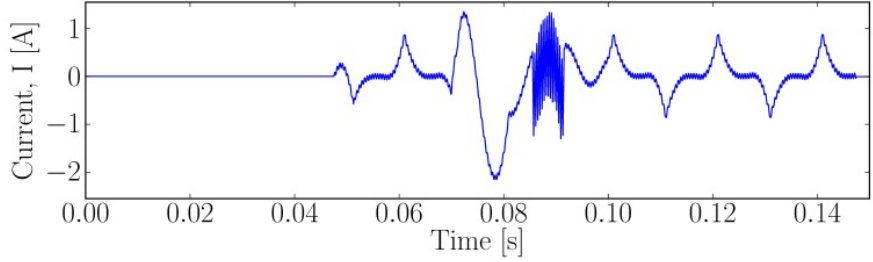

(a)

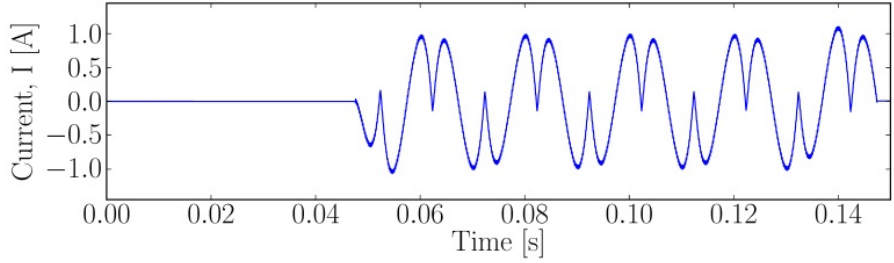

(b)

Figure 3. Residual currents detected by EMD after the APF is turned on: (a) residual currents at line 1; and (b) residual currents at line 3.

The residual current components detected by EMD for line 1 and line 3, are shown in Figure 3.

To testify the effectiveness of the proposed APF control, the results obtained is compared with the ones obtained with the conventional $p$-$q$ theory. The results obtained are shown in Figure 4 and 5, respectively.

- **Source current waveforms:** from Figure 4a and 5a, both methods are able to produce sinusoidal currents after the APF is turned on. However, the current disturbance which was injected into line 1 is still present for the APF which employed the modified $p$-$q$ theory. On the other hand, this disturbance is not observed in Figure 4a.
- **Active and reactive powers:** both methods are able to give constant active power and compensate the reactive power after the APF is turned on. However, for the case which the APF uses only the modified $p$-$q$ theory, there is a disturbance in both the active and reactive powers at the duration in which the current disturbance is introduced. The reactive power rises to 100 W in that period.

- **Power factor:** the proposed method is able to achieve unity power factor. The oscillation caused by the non-linear non-stationary current is totally removed. With the conventional modified *p-q* theory, the power factor reduces to 0.85, due to the algorithm is not able to eliminate the non-linear non-stationary current disturbance.
- **Neutral current:** both methods are able to give zero neutral current after the APF is turned on.

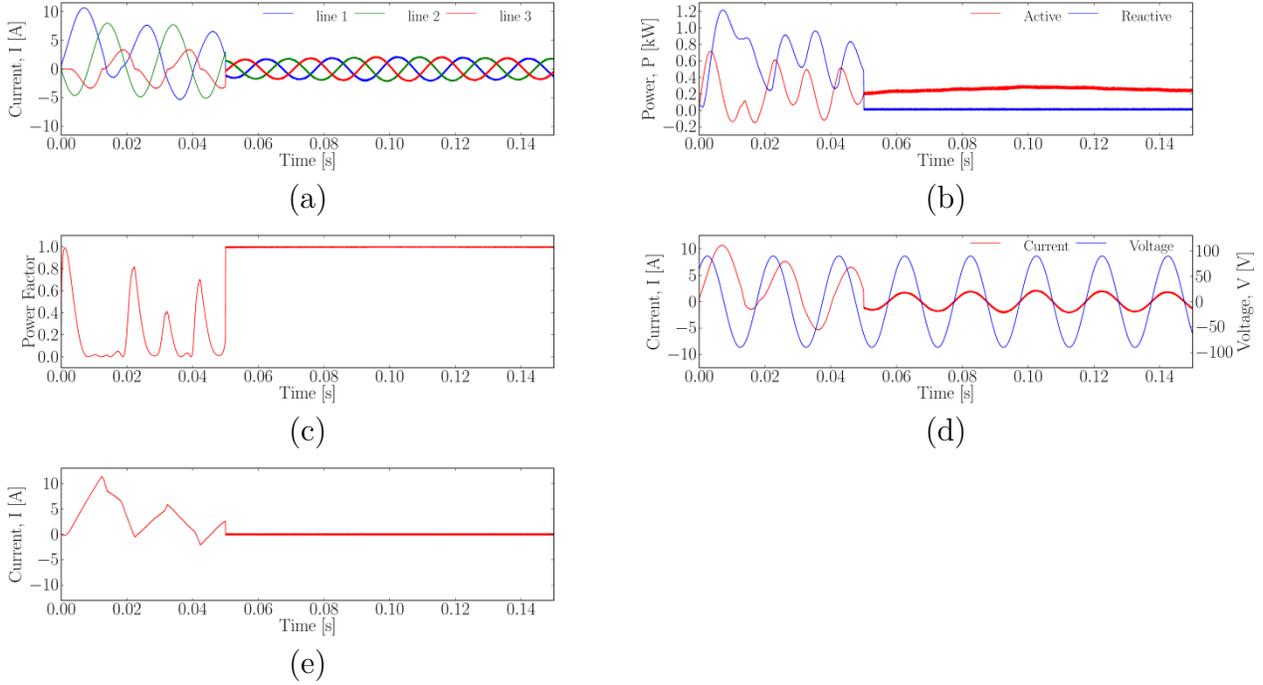

Figure 4. Compensation results with the proposed method: (a) source currents; (b) source active and reactive powers; (c) power factor; (d) line 1 source current and voltage; and (e) neutral current.

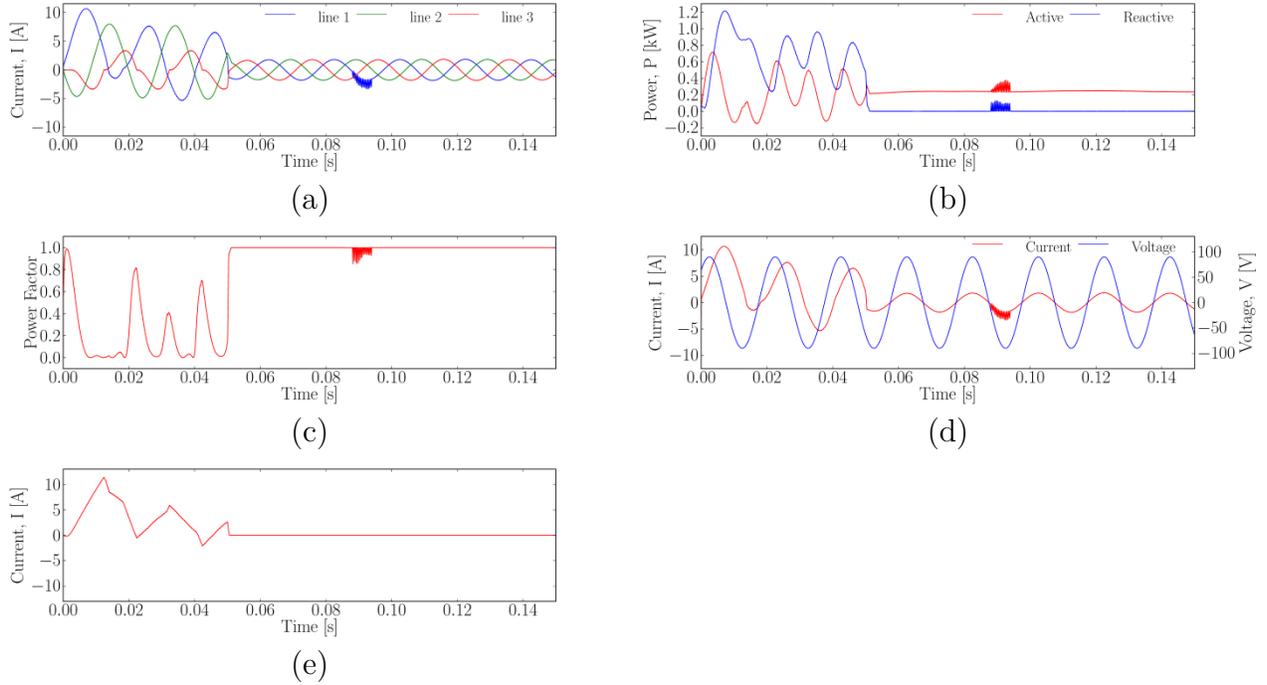

Figure 5. Compensation results with modified *p*-*q* theory: (a) source currents; (b) source active and reactive powers; (c) power factor; (d) line 1 source current and voltage; and (e) neutral current.

The results confirmed the conventional *p*-*q* theory is not able to compensate non-linear non-stationary currents. The effects of non-linear non-stationary current disturbance are not observed with the proposed method.

## 5  Conclusions

In the presence of non-periodic currents, [Watanabe et al., 2000] showed that *p*-*q* theory cannot give constant power and sinusoidal current at the same time. However, with proper formulation and signal decomposition technique, modified *p*-*q* theory can give sinusoidal source current and unity power factor, even with the presence of non-linear non-stationary currents. EMD is effective in singulating the residual currents from the load currents. The definition of current components used in this study to formulate the APF control fits the current definition given by [Czarnecki, 2000]. Further investigation is needed to cover the compensation which involved voltage distortions.

EMD is a promising technique for non-stationary and non-linear signal processing due to its adaptive nature [Huang et al., 1998]. This feature makes it suitable to be used in the control applications for power systems such as in the field of wide-area monitoring and control.

## Acknowledgements

The authors wish to thank the Norwegian University of Science and Technology for all the support given in this research.

Togasawa S., Murase T., Nakano H. and Nabae A. Reactive power compensation based on a novel cross-vector theory. Trasactions on Industrial Application of the IEEE Japan 114 (3): 340-341, 1994.

Tlusty J., Svec J., Sendra J. B. and Valouch V. Analysis of generalized non-active power theory for compensation of non-periodic disturbances. In Proceedings of the International Conference on Renewable Energies and Power Quality. Santiago de Compostela, Spain, March 28-30, 2012.

Tolbert L. M., Xu Y., Chen J., Peng F. Z. and Chiasson J. N. Application of compensators for non-periodic currents. IEEE Power Electronic Letters 1 (2): 45-50, 2003.

Qi L., Qian L., Cartes D. and Woodruff S. Initial results in Prony analysis for harmonic selective active filters. In Proceedings of the Power Engineering Society General Meeting. Montreal, Canada, June 18-22, 2006.

Watanabe E. H. And Aredes M. Compensation of non-periodic currents using the instantaneous power theory. In Proceedings of the 2000 IEEE Power Engineering Society Summer Meeting, pp. 994-999. Seattle, USA, July 16-20, 2000.12